\documentstyle[twoside,psfig,amsgen,amsmath,amssymb,curvesls,epic,graphicx,psfrag]{article}

\oddsidemargin=\evensidemargin
\newcommand{\C}{{\bf{C}}}
\newcommand{\I}{\bf{I}}
\newcommand{\Real}{{\bf{R}}}
\newcommand{\Half}{{\bf{H}}}
\newcommand{\Zed}{{\bf{Z}}}
\newcommand{\AT}{{\bf{A}}}

\newcommand{\Lau}{{\Zed[q^{\pm 1},t^{\pm 1}]}}
\newcommand{\imes}{{\hskip -.6mm\times\hskip -.6mm}}

\newcommand{\GL}{{\mathrm{GL}}}
\newcommand{\Hom}{\mathrm{Hom}}
\newcommand{\Diff}{\mathrm{Diff}}

\newcommand{\Ho}{{\mathrm{H}}}

\newcommand{\F}{{\mathcal{F}}}
\newcommand{\Config}{{\mathcal{C}}}
\newcommand{\La}{{\mathcal {L}}}
\newcommand{\PDisc}{{\mathcal{P}}}

\catcode`\@=11
\long\def\@makefntext#1{
\protect\noindent \hbox to 3.2pt {\hskip-.9pt
$^{{\eightrm\@thefnmark}}$\hfil}#1\hfill}		

\def\ps@myheadings{%
    \let\@oddfoot\@empty\let\@evenfoot\@empty
    \def\@evenhead{\slshape\leftmark\hfil}
    \def\@oddhead{\hfil{\slshape\rightmark}}
    \let\@mkboth\@gobbletwo
    \let\sectionmark\@gobble
    \let\subsectionmark\@gobble
    }

\catcode`@=11		      		     
\def\ps@plain{\let\@mkboth\@gobbletwo
     \def\@oddhead{}\def\@oddfoot{\eightrm\hfil\thepage
     \hfil}\def\@evenhead{}\let\@evenfoot\@oddfoot}

\renewcommand{\thefootnote}{\fnsymbol{footnote}}

\newcounter{sectionc}\newcounter{subsectionc}\newcounter{subsubsectionc}
\renewcommand{\section}[1] {\vspace{12pt}\addtocounter{sectionc}{1}
\setcounter{subsectionc}{0}\setcounter{subsubsectionc}{0}\noindent
	{\tenbf\thesectionc. #1}\par\vspace{5pt}}
\renewcommand{\subsection}[1] {\vspace{12pt}\addtocounter{subsectionc}{1}
	\setcounter{subsubsectionc}{0}\noindent
	{\bf\thesectionc.\thesubsectionc.
	{\kern1pt \bfit #1}}\par\vspace{5pt}}
\renewcommand{\subsubsection}[1] {\vspace{12pt}
	\addtocounter{subsubsectionc}{1}
	\noindent
	{\tenrm\thesectionc.\thesubsectionc.\thesubsubsectionc.	{\kern1pt
	\it #1}}\par\vspace{5pt}}

\newcounter{appendixc}
\newcounter{subappendixc}[appendixc]
\newcounter{subsubappendixc}[subappendixc]

\renewcommand{\appendix}[1] {\vspace{12pt}	
	\refstepcounter{appendixc}		
	\setcounter{figure}{0}
	\setcounter{table}{0}
	\setcounter{lemma}{0}
	\setcounter{theorem}{0}
	\setcounter{corollary}{0}
	\setcounter{definition}{0}
	\setcounter{equation}{0}
	\renewcommand{\thefigure}{\Alph{appendixc}.\arabic{figure}}
	\renewcommand{\thetable}{\Alph{appendixc}.\arabic{table}}
	\renewcommand{\theappendixc}{\Alph{appendixc}}
	\renewcommand{\thelemma}{\Alph{appendixc}.\arabic{lemma}}
	\renewcommand{\thetheorem}{\Alph{appendixc}.\arabic{theorem}}
	\renewcommand{\thedefinition}{\Alph{appendixc}.\arabic{definition}}
	\renewcommand{\thecorollary}{\Alph{appendixc}.\arabic{corollary}}
	\renewcommand{\theequation}{\Alph{appendixc}.\arabic{equation}}
	\noindent{\tenbf Appendix \theappendixc #1}\par\vspace{5pt}}

\newcommand{\textlineskip}{\baselineskip=13pt}
\newcommand{\smalllineskip}{\baselineskip=10pt}

\newcommand{\copyrightheading}[1]
	{\vspace*{-2.5cm}\smalllineskip{\flushleft
	{\footnotesize Journal of Knot Theory and Its Ramifications #1}\\
   	{\footnotesize \copyright\kern2pt World Scientific
         Publishing Company}\\
         }}

\newcommand{\publisher}[2]{{\begin{center}\footnotesize\smalllineskip
	Received #1\\
	Revised #2
        \end{center}
	}}

\def\abstracts#1#2#3#4{{
	\centering{\begin{minipage}{4.5in}\footnotesize\baselineskip=10pt
	\centerline{ABSTRACT}
	\parindent=15pt #1\par
	\parindent=15pt #2\par
	\parindent=15pt #3\par
	\parindent=15pt #4\par
	\end{minipage}}\par}}

\renewenvironment{thebibliography}[1]
	{\frenchspacing
	 \ninerm\baselineskip=11pt
	 \begin{list}{[\arabic{enumi}]}
	{\usecounter{enumi}\setlength{\parsep}{0pt}
	 \setlength{\leftmargin 13.7pt}{\rightmargin 0pt} 
	 \setlength{\itemsep}{0pt} \settowidth
	{\labelwidth}{[#1]}\sloppy}}{\end{list}}

\newcounter{itemlistc}
\newcounter{romanlistc}
\newcounter{alphlistc}
\newcounter{arabiclistc}

\newcommand{\fcaption}[1]{
        \refstepcounter{figure}
        \setbox\@tempboxa = \hbox{\footnotesize Fig.~\thefigure. #1}
        \ifdim \wd\@tempboxa > 5in
           {\begin{center}
        \parbox{5in}{\footnotesize\smalllineskip Fig.~\thefigure. #1}
            \end{center}}
        \else
             {\begin{center}
             {\footnotesize Fig.~\thefigure. #1}
              \end{center}}
        \fi}

\newcommand{\tcaption}[1]{
        \refstepcounter{table}
        \setbox\@tempboxa = \hbox{\footnotesize Table~\thetable. #1}
        \ifdim \wd\@tempboxa > 5in
           {\begin{center}
        \parbox{5in}{\footnotesize\smalllineskip Table~\thetable. #1}
            \end{center}}
        \else
             {\begin{center}
             {\footnotesize Table~\thetable. #1}
              \end{center}}
        \fi}

\def\pmb#1{\setbox0=\hbox{#1}
	\kern-.025em\copy0\kern-\wd0
	\kern.05em\copy0\kern-\wd0
	\kern-.025em\raise.0433em\box0}

\def\fnt#1#2{\footnotetext{\kern-.3em
	{$^{\mbox{\scriptsize #1}}$}{#2}}}

\font\tenrm=cmr10
 
\font\tenbf=cmbx10
\font\bfit=cmbxti10 at 10pt
\font\ninerm=cmr9

\font\eightrm=cmr8

\def\@begintheorem#1#2{\trivlist	
	\item[\hskip\labelsep{\bf #1\ #2.}]}
\def\@opargbegintheorem#1#2#3{\trivlist
	\item[\hskip\labelsep{\bf #1\ #2\ (#3).}]}

    	{\setcounter{itemlistc}{0}		
	 \begin{list}{$\bullet$}		
	{\usecounter{itemlistc}			
	 \leftmargin10pt	       
	 \setlength{\parsep}{0pt}
	 \setlength{\itemsep}{0pt}     
	}}{\end{list}}

	{\setcounter{romanlistc}{0}		
	 \begin{list}{$($\roman{romanlistc}$)$}	
	{\usecounter{romanlistc}		
	 \leftmargin18pt 
	 \setlength{\parsep}{0pt}
	 \setlength{\itemsep}{0pt}	
	 \settowidth{\labelwidth}{#1}                          
	}}{\end{list}}

	{\setcounter{enumii}{0}			
	 \begin{list}{$($\alph{enumii}$)$}	
	{\usecounter{enumii}			
	 \leftmargin18pt		
	 \setlength{\parsep}{0pt}
	 \setlength{\itemsep}{0pt}	
	 \settowidth{\labelwidth}{#1}
	}}{\end{list}}

\def\qed{\hbox{${\vcenter{\vbox{			
   \hrule height 0.4pt\hbox{\vrule width 0.4pt height 6pt
   \kern5pt\vrule width 0.4pt}\hrule height 0.4pt}}}$}}

\renewcommand{\thefootnote}{\fnsymbol{footnote}}  

\def\theequation{\thesectionc.\arabic{equation}}  

\begin{document}
\setlength{\textheight}{7.7truein}  

\thispagestyle{empty}

\markboth{\protect{\footnotesize\it On the Image of the Lawrence-Krammer representation}}{\protect{
           \footnotesize\it On the Image of the Lawrence-Krammer representation}}

\normalsize\textlineskip

\setcounter{page}{1}

\copyrightheading{}	

\vspace*{0.88truein}

\centerline{\bf ON THE IMAGE OF THE}
\baselineskip=13pt
\centerline{\bf LAWRENCE-KRAMMER REPRESENTATION}
\vspace*{0.37truein}
\centerline{\footnotesize RYAN D. BUDNEY}
\baselineskip=12pt
\centerline{\footnotesize\it Max Planck Institute for Mathematics}
\baselineskip=10pt
\centerline{\footnotesize\it Vivatsgasse 7.}
\baselineskip=10pt
\centerline{\footnotesize\it D-53111 Bonn}
\baselineskip=10pt
\centerline{\footnotesize\it Germany}

\vspace*{0.225truein}
\publisher{(Leave 1 inch blank space for publisher.)}

\vspace*{0.21truein}
\abstracts{A non-singular sesquilinear form is constructed
that is preserved by the Lawr\-ence-Kram\-mer representation.
It is shown that if the polynomial variables $q$ and $t$ of the
Lawr\-ence-Kram\-mer representation are chosen to be appropriate
algebraically independent unit complex numbers,
then the form is negative-definite Hermitian.  Using
the fact that non-invertible knots exist this result implies
that there are matrices in the image of the Lawrence-Krammer
representation that are conjugate in the unitary group, yet the braids
that they correspond to are not conjugate as braids. The two primary
tools involved in constructing the sesquilinear form are Bigelow's
interpretation of the Lawrence-Krammer representation, together with
the Morse theory of functions on manifolds with corners.}{}{}{}



\vspace*{1pt}\textlineskip	
\section{Introduction}	
\vspace*{-0.5pt}

This paper takes a Morse-theoretic approach to the Lawrence-Krammer
representation.  The Lawrence-Krammer representation is an injective homomorphism
$B_n \to \GL_{n \choose 2} \Lau$ \cite{Bi1} \cite{Kr}. The representation
has a natural description as the action
of the braid group on the middle-dimensional homology of a certain four-dimensional
manifold, where the homology is thought of as a module over a Laurent polynomial
ring, which enters the picture as the group ring of a free abelian group of
covering transformations similar to the Burau representation.
The middle dimensional homology
of any even-dimensional manifold has an intersection product pairing, and
this is used to
construct a sesquilinear form that is preserved by the Lawrence-Krammer
representation, analogously to the work of Long \cite{Long} and similar work for the Gassner
and Burau representations by Abdulrahim \cite{Ab} and Squier \cite{Sq}.
Using the Morse theory of functions on
manifolds with corners developed by Handron \cite{Ha1} \cite{Ha2}, this sesquilinear
pairing is explicitly computed in Section 4.

The embeddings of $\Lau \to \C$ are parametrized by algebraically
independant $q,t \in \C$,
so we can think of the Lawrence-Krammer representation as a map
$\C^2 \to \Hom(B_n, GL_{n \choose 2} \C)$.
Provided $|q|=|t|=1$ the sesquilinear form that is preserved
by the Lawrence-Krammer representation is Hermitian. It
is shown in Theorem 2 that for certain values of
$t$ and $q$ this Hermitian
form is negative definite, and thus the image of the
Lawrence-Krammer representation
$B_n \to GL_{n \choose 2} \C$ has compact closure.  Using this
result, we address the question of conjugacy in the image of the Lawrence-Krammer
representation in Section 5. This result may be of
interest to Braid Cryptographers
\cite{An} as it gives insight into the difference between the conjugacy
problem in braid groups versus the conjugacy problem in the target
matrix group, $GL_{n \choose 2} \C$.

\section{Generalities on the Lawrence-Krammer representation}

This section begins with the definition of the Lawrence-Krammer
representation and more generally the Lawrence representations
and the sesquilinear forms that they preserve.
\vskip 7mm

{\bf Definition~1.} The {\bf configuration space} $\Config_n X$ of
$n$ points in a topological space $X$ is the space
$(X^n - \Delta_nX) / S_n$. Here
$\Delta_nX = \{ (x_1,\cdots,x_n) \in X^n : x_i=x_j \text{ for some } i\neq j\}$.
$S_n$ is the symmetric group, acting by permuting the factors of the
product $X^n$.

{\it A convention on the fundamental group
that is used throughout this paper is that if $f, g : \I=[0,1] \to X$
are loops, then the concatenation $fg$ denotes the
loop such that $(fg)(t)=g(2t)$ for $0 \leq t \leq \frac{1}{2}$
and $(fg)(t)=f(2t-1)$ for $\frac{1}{2} \leq t \leq 1$.}

The {\bf braid group}
$B_n := \pi_0 \Diff (D^2,n)$ is the
mapping class group of a disc with $n$ marked points in the interior,
where the diffeomorphisms restrict to the identity on the boundary.
An equivalent definition of the braid group is the
fundamental group of the configuration space of $n$ points in a
disc, $B_n = \pi_1 \Config_n D^2$.
The fact that the two definitions are equivalent is an easy
consequence of the homotopy long exact sequence of the fibration
$\Diff (D^2) \to \Config_n D^2$, together with Smale's theorem that
$\Diff (D^2)$ is contractible \cite{Sm}. The map
$\Diff (D^2) \to \Config_n D^2$
is given by fixing a configuration in $\Config_n D^2$ and evaluating
it on a diffeomorphisms of $\Diff(D^2)$.
The fact that evaluation maps are fibrations is due to Palais \cite{Pa}.
The boundary map $\pi_1 \Config_n D^2 \to \pi_0 \Diff (D^2,n)$ is
a homomorphism with the above concatenation convention in 
$\pi_1 \Config_n D^2$.

$\PDisc_n$ will denote the closed unit disc with $n$ interior points
removed. Let
$ab : B_i \to \Zed$ for $i \in \{1, 2, 3, \cdots \}$ be the
abelianization maps. A convention in this paper is that
$ab(\sigma_i)=1$, where $\sigma_i$ is the half Dehn twist in Figure 1.
Let $T :  \pi_1 \Config_k \PDisc_n \to \Zed$ be the composite of
the forgetful map $\pi_1 \Config_k \PDisc_n \to \pi_1 \Config_k D^2$ with
the abelianization map, and similarly let $R : \pi_1 \Config_k \PDisc_n \to \Zed$
be the composite of the inclusion map
$\pi_1 \Config_k \PDisc_n \to \pi_1 \Config_{k+n} D^2$
with the abelianization map, and define $Q : \pi_1 \Config_k \PDisc_n \to \Zed$ by
the identity $Q(f):={\frac{R(f)-T(f)}{2}}$. Let $\La \Config_k \PDisc_n$
be the abelian Galois covering space
of $\Config_k \PDisc_n$ such that the image of the map
$\pi_1 \La \Config_k \PDisc_n \to \pi_1 \Config_k \PDisc_n$
is $ker(Q)\cap ker(T)$. The $k$-th Lawrence representation of $B_n$, as
described by Bigelow in \cite{Bi3} is the action of
$B_n$ on $\Ho_k (\La \Config_k \PDisc_n)$.

\begin{center}Figure 1
{
\psfrag{s}{$\sigma_i$}
\psfrag{bs}{$\partial^{-1}\sigma_i$}
\psfrag{1}{\scalebox{0.7}{$1$}}
\psfrag{i}{\scalebox{0.7}{$i$}}
\psfrag{i+1}{\scalebox{0.7}{$i+1$}}
\psfrag{n}{\scalebox{0.7}{$n$}}
$$\includegraphics[width=10cm]{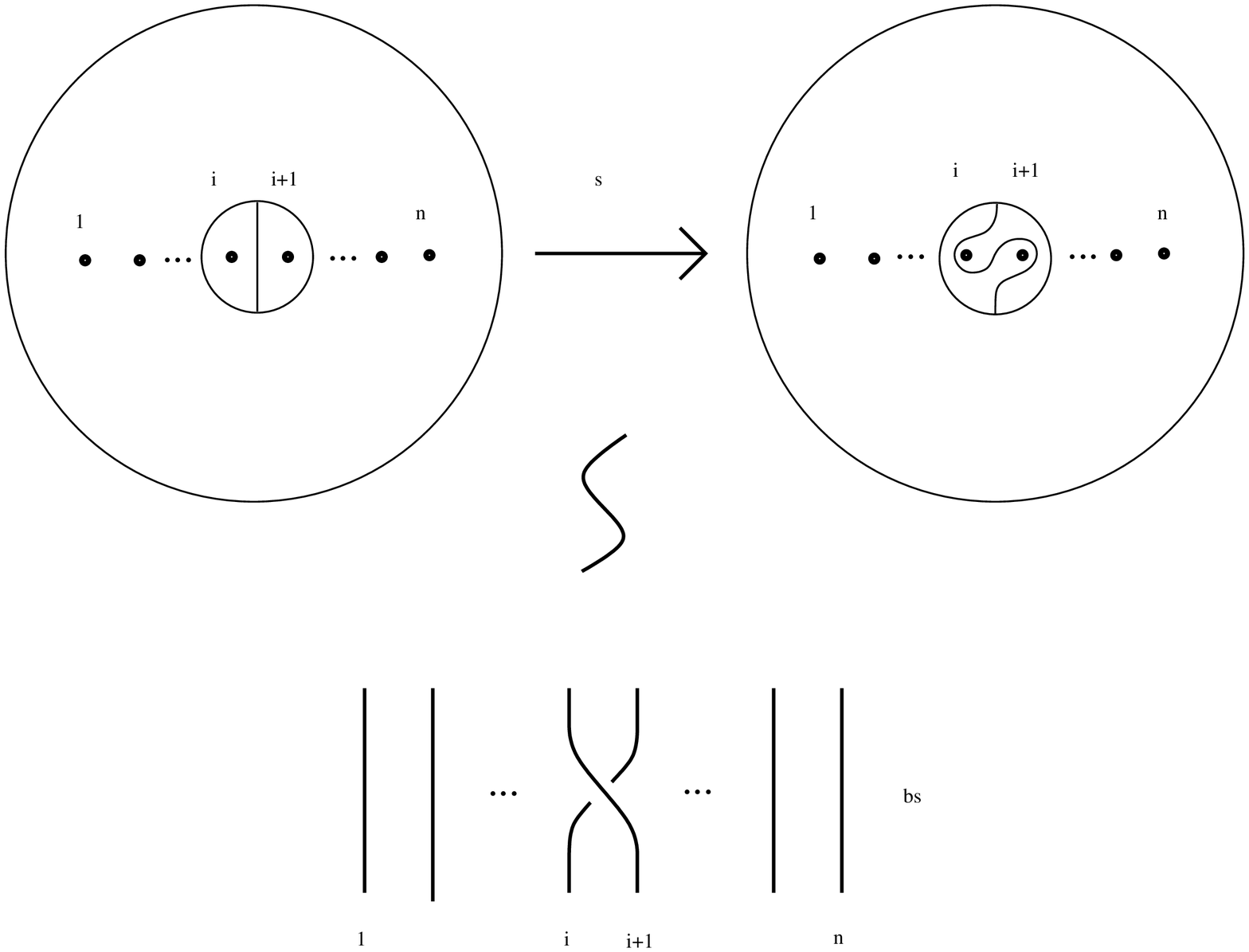} $$
}
An example of the isomorphism $\partial : \pi_1 \Config_n D^2 \to \pi_0 \Diff(D^2,n)$
\end{center}
\vskip 7mm

$\La \Config_k \PDisc_n \to \Config_k \PDisc_n$
is a `normal' or Galois cover.  
For $k \geq 2$ its group of covering transformations is precisely 
$\Zed \times \Zed$ and
can be identified with the image of $Q \times T$.
$q$ and $t$ will denote
the covering transformations corresponding to $1 \times 0$ and $0 \times 1$
respectively in the image of $Q \times T$.
With these definitions, $\Ho_2 (\La \Config_2 \PDisc_n)$ is a module over
the Laurent polynomial ring $\Lau$, which is the group ring of the
group of covering transformations $\langle q, t \rangle = \Zed\times\Zed$.

In \cite{Bi1}, the Lawrence-Krammer representation was defined
as the action of the braid group $B_n$ on
$\Ho_2 (\La \Config_2 \PDisc_n) \otimes_{\Lau} \F$
where $\F$ is some field containing $\Lau$ such as $\C$.
It has since been shown that $\Ho_2 (\La \Config_2 \PDisc_n)$
is free over the Laurent polynomial ring $\Lau$ \cite{Pl}. Unfortunately, those
free generators are not convenient to use, so in this paper we
will restrict the Lawrence-Krammer representation to
a certain full-rank, free, $B_n$-invariant submodule of
$\Ho_2 (\La \Config_2 \PDisc_n)$ which will be described precisely
in Section 3.

{\bf Definition~2.}
The {\bf total intersection product}
is a sesquilinear form
$$\langle \cdot,\cdot \rangle : \Ho_2 (\La \Config_2 \PDisc_n) \oplus \Ho_2 (\La \Config_2 \PDisc_n) \to \Lau$$
defined by
$$\langle v,w \rangle := \sum_{(i,j) \in \Zed \imes \Zed} \mu(v,t^iq^jw)t^iq^j$$
where $\mu : \Ho_2 (\La \Config_2 \PDisc_n) \oplus \Ho_2 (\La \Config_2 \PDisc_n) \to \Zed$ is the intersection
product. See \cite{ST} or \cite{Br} for definitions of the
intersection product on manifolds, and \cite{Lan} for basics on sesquilinear
forms. Sometimes this pairing is called a Blanchfield form \cite{Kaw}.

%
To compute the intersection product, one could take the CW-decomposition
given in Bigelow's paper \cite{Bi1} and notice that all homology classes are realizable
by compact surfaces. Unfortunately, this is potentially very difficult,
as Bigelow's generators $v_{i,j}$ are genus two surfaces and rather
difficult to visualize, moreover they are not transverse.
To bypass this difficulty, we compute the intersection product using two easy to
visualize transverse CW-decompositions of $\Config_2 \PDisc_n$ that
come from a Morse function on $\Config_2 \PDisc_n$.

\section{A little Morse theory}
In section 3.1 we review Morse theory on manifolds with corners.  In
section 3.2 the Morse theory is applied to get dual CW-decompositions
of the configuration space of two points in a planar surface, and in
3.3 we apply these results to study the Lawrence-Krammer representation.

\subsection{A survey of Morse theory on manifolds with corners}

Morse theory on manifolds with corners
has been studied for some time, although it is not a commonly known
branch of Morse theory.  A good general reference for Morse theory
is Milnor's book \cite{Mi}, and for Morse theory on
manifolds with corners, Handron's papers
\cite{Ha1} \cite{Ha2}. A summary of the relevant
theory is given below.

{\bf Definition~3.} A smooth
$n$-dimensional {\bf manifold with corners} is a $2^{\text{nd}}$-countable
Hausdorff topological space $X$ together with a family
of maps $\AT$ where if $\xi \in \AT$ then $\xi : U \to \Half^n_j$
is a homeomorphism between an open subset $U$ of $X$ and $\Half^n_j$
for some $j \in \{0,1,\cdots,n\}$.
$$\Half^n_j = \{ w \in \Real^n : w \cdot e_i \geq 0 \ \forall \ 1\leq i \leq j \}$$
$\xi \in \AT$ is called a {\bf chart}.
We demand that the union of the domains of the charts in $\AT$ is $X$, and
if any two charts
$\xi,\psi$ have overlapping domains $\xi : U \to \Half^n_j$ and
$\psi : V \to \Half^n_k$
then they must be smoothly compatible in the sense that
 $\xi \circ \psi^{-1}_{|\psi (U\cap V)}$ must be a smooth diffeomorphism
 from $\psi (U \cap V)$ to $\xi (U \cap V)$ \cite{GP}. This allows
 us to define smooth functions between manifolds with corners and
 derivatives of such functions analogously to \cite{Hirsch}.

The {\bf $i$-dimensional strata} $X^{(i)}$ of $X$ is the set of points
$x \in X$ such that there exists a chart $\xi : U \to \Half^n_{n-i}$, $x \in U$
with $\xi(x) \cdot e_j = 0 \ \forall \ 1 \leq j \leq n-i$.
A {\bf critical point} of a function $f : X \to \Real$ is a point
$x \in X$ such that if $x \in X^{(i)}$ then $Df_{|X^{(i)}}(x)=0$.
A critical point $x$ of $f$ is {\bf non-degenerate} if the Hessian
matrix $D^2f_{|X^{(i)}}(x)$ is non-singular and if for all $v \in T_xX$ that
point into the strata $X^{(i+1)}$, $Df_x(v)\neq 0$. A function $f : X \to \Real$
is a {\bf Morse function} if all of its critical points are non-degenerate.
A non-degenerate critical point is {\bf essential} if for all $v \in T_xX$ that
point into $X^{(i+1)}$,  $Df_x(v) > 0$.

{\bf Theorem~1.}
\cite{GM} \cite{V} \cite{Ha1} \cite{Ha2} Given a Morse function
$f : X \to \Real$ the homotopy type of $f^{-1}(-\infty,c]$ changes
only at $c=f(x)$ for $x$ an essential
critical point of $f$. Provided $f^{-1}(c)$ contains only one essential
critical point, $f^{-1}[c-\epsilon,c+\epsilon]$ is homotopy
equivalent to $f^{-1}(c-\epsilon)$ union a cell, the dimension of the
cell is given by the index of $D^2f_{|X^{(i)}}(x)$.

\begin{center}
\begin{minipage}{10cm}
\begin{center} Figure 2
{
\psfrag{f}{$f$}
\psfrag{R}{$\Real$}
\psfrag{X1}{$X^2$}
\psfrag{X2}{$X^1$}
\psfrag{x}{$x$}
\psfrag{v}{$v$}
\begin{center}
$$\includegraphics[width=10cm]{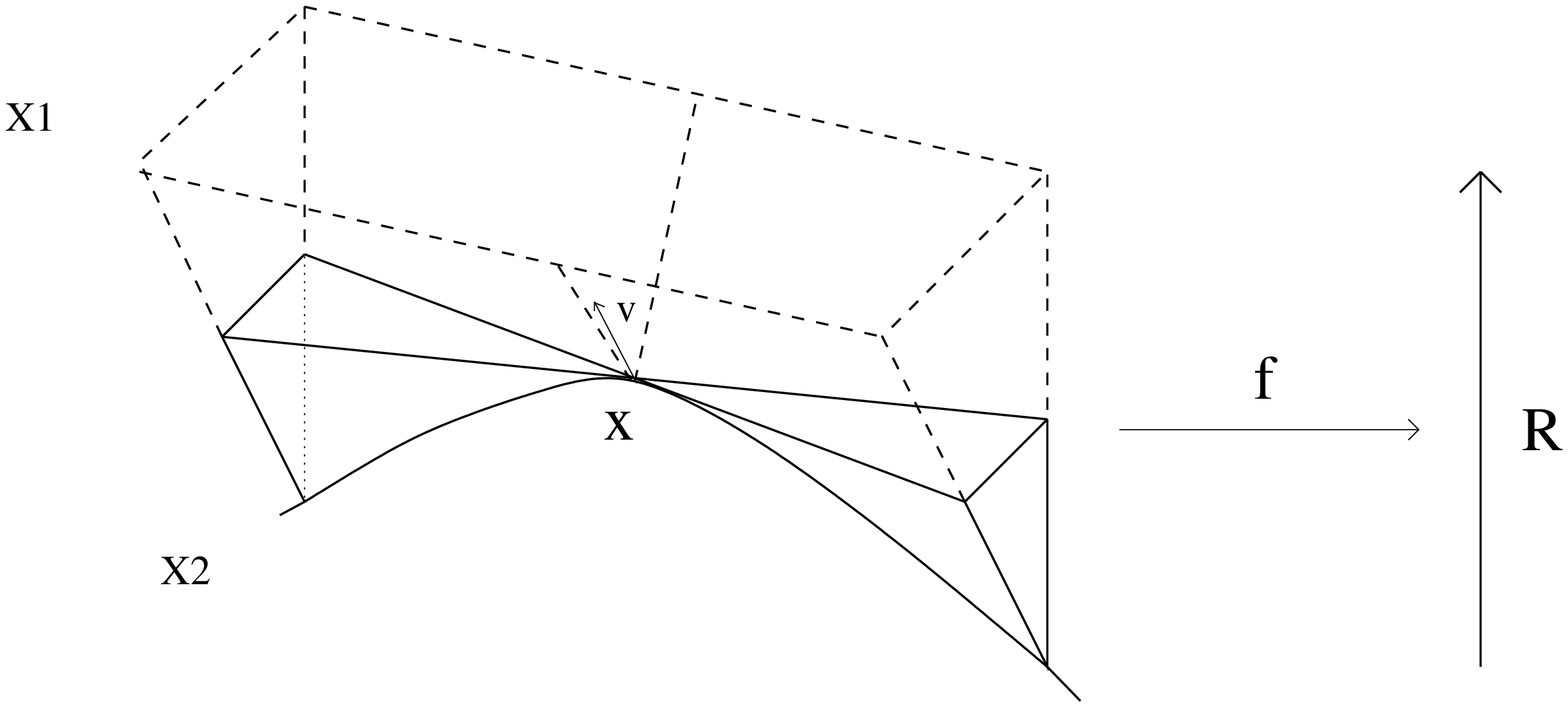} $$
\end{center}
Essential critical point on a $1$-dimensional strata resulting in a
$1$-cell attachment.
}
\end{center}
\end{minipage}
\end{center}

Given $M$ a smooth manifold with boundary, $M^n$ is naturally a smooth
manifold with corners. Thus, $\Config_n M$ is a smooth manifold with
corners.

{\bf Lemma~1.}
Let $R$ be a sub-manifold of $\Real^2$, and consider the
distance function $d : \Config_2 R \to \Real$, $d([(z_1,z_2)])=|z_1-z_2|$
where $| \cdot |$ is the standard Euclidean norm on $\Real^2$.
A point
$[z_1,z_2] \in \Config_2 R$ is a critical point of $d$ if and only
if $[z_1,z_2]$ sits on the 2-dimensional strata of $\Config_2 R$ and
$z_1 - z_2$ is perpendicular to both $T_{z_1} \partial R$ and
$T_{z_2} \partial R$.  The Hessian is non-singular at $[z_1,z_2]$ provided
the polynomial $xy+x+y$ does not have a root at
$(x,y)=(\langle k_1,z_1-z_2\rangle,\langle z_2-z_1,k_2\rangle)$
where $k_1$ and $k_2$ are the curvature vectors of $\partial R$ at
$z_1$ and $z_2$ respectively. The index of the Hessian is the number
of roots $t\in (0,\infty)$ of the polynomial
$xy+x+y$ where $(x,y)=(\langle k_1,z_1-z_2\rangle+t,\langle z_2-z_1,k_2\rangle+t)$.

{\bf Proof.} $d$ has critical points only on the 2-dimensional strata of
$\PDisc_n$ since if $[z_1,z_2]$ is not on the $2$-dimensional strata,
at least one of $z_1$ or $z_2$ must be in the interior of $R$ and so
the derivative of $d$ must be non-zero.
Let $[z_1,z_2] \in (\Config_2 R)^{(2)}$, thus both $z_1$ and $z_2 \in \partial R$.
Let $f : (-\epsilon,\epsilon) \to \partial R$
and $g : (-\epsilon, \epsilon) \to \partial R$ be arclength-preserving parametrizations
of neighborhoods of $z_1$ and $z_2$ in $\partial R$ respectively with
$f(0)=z_1$ $g(0)=z_2$.
The critical points of $d \circ [f(x),g(y)]$ are precisely the same as
the critical points of the square of the distance function, which is
a polynomial function
$d^2 \circ [f(x),g(y)] = \langle f(x)-g(y),f(x)-g(y) \rangle$.
A quick computation shows that the derivative of the above polynomimal is the
$1\times 2$ matrix $[\langle f'(x),f(x)-g(y) \rangle, \langle g(y)-f(x),g'(y) \rangle]$.
To prove the statement about non-singularity and index of the Hessian, we compute
$D^2(d^2\circ[f(x),g(y)])_{(0,0)}$.
This is the matrix
$\left[
\begin{matrix}
 c_1 +1 & -1 \\ -1 & c_2+1
 \end{matrix}\right]$
 where $c_1 = \langle f''(0),f(0)-g(0)\rangle$ and
  $c_2 = \langle g(0)-f(0),g''(0)\rangle$. The result follows.
\qed
\begin{center}
\begin{minipage}{10cm}
\begin{center} Figure 3
{
\psfrag{x1}{$y=-1$}
\psfrag{y1}{$x=-1$}
\psfrag{x}{$x$}
\psfrag{y}{$y$}
\begin{center}
$$\includegraphics[width=7cm]{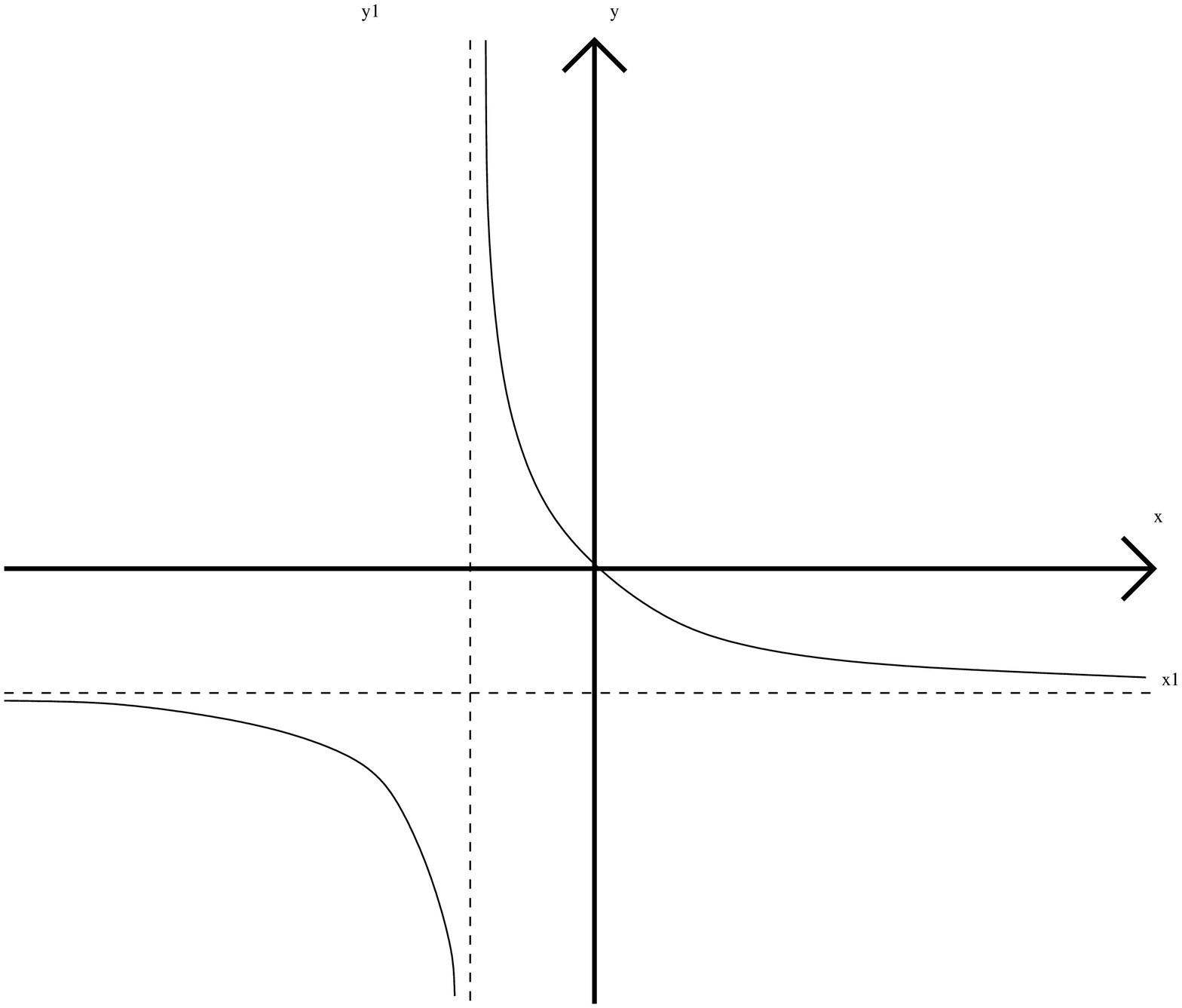} $$
\end{center}
Roots of the polynomial $xy+x+y$
}
\end{center}
\end{minipage}
\end{center}

\subsection{A CW-decomposition}
For the purposes of the remainder of this paper, consider the punctured disc
$\PDisc_n$ to be the genus zero compact connected 2-manifold with $n+1$ boundary
components.  Consider an embedding of $\PDisc_n$ in $\Real^2$ where the boundary
consists entirely of ellipses, as in Figure 4.

{\bf Lemma~2.}
\label{cw1}
The distance
function $d : \Config_2 \PDisc_n \to \Real$ is Morse.
The CW decomposition given by the Morse function $d$ has one $0$-cell,
$2n+1$ $1$-cells, labeled $x_i, y_i, i \in \{ 1, \cdots , n \}$ and $b$.
There are also ${n\choose 2} + 2n$ $2$-cells, with attaching maps:

$$
\begin{array}{llr}
 Z_{i,j} & [bx_ib^{-1},x_j^{-1}] & 1\leq i<j\leq n \\
 Z_i & y_i^{-1}b^{-1}x_ibx_i & 1\leq i\leq n \\
 Z'_i & [b,y_i] & 1\leq i\leq n
\end{array}
$$
We use the convention that $[a,b]=aba^{-1}b^{-1}$.

{\bf Proof.}
$d^{-1}(-\infty,x]$ for $x > 0$ small is has
the homotopy type of the projectivised unit tangent bundle of $\PDisc_n$,
which is diffeomorphic to $S^1 \times \vee_{i=1}^{n} S^1$.
We give $S^1 \times \vee_{i=1}^n S^1$ the product CW-decomposition,
with counter-clockwise orientations given to the 1-cells.
The cell corresponding to the
$S^1$ factor will be denoted $b$ and
the $n$ wedge summands will be denoted $y_i$ for $i \in \{ 1, \cdots, n \}$.
The product structure gives us the $2$-cells $Z_i'$ with attaching maps
$[b,y_i]$.

There are $4{n \choose 2} + 10n + 2$ critical points, but
only ${n \choose 2} + 2n$ of them are essential.  If one
restricts $d$ to each connected 2-dimensional stratum,
the essential critical points are the maxima on each stratum, and
if both points of the configuration lay on the same ellipse, then
there is an additional essential critical point corresponding to the
minor axis of the ellipse.

Of the ${n \choose 2} + 2n$ essential critical points, ${n \choose 2} + n$
are $2$-cell attachments (the maxima of $d$ on each connected 2-dimensional stratum),
and the remaining $n$ essential critical points are $1$-cell attachments.
This is true because of Lemma 1 tells us that a critical point corresponding to
the minor axis of an ellipse is a saddle point.

\begin{center}
\begin{minipage}{10cm}
\begin{center} Figure 4.
{
\psfrag{ZIJ}{$Z_{i,j}$'s}
\psfrag{zi}{$Z_i$'s}
\psfrag{xis}{$x_i$'s}
$$\includegraphics[width=10cm]{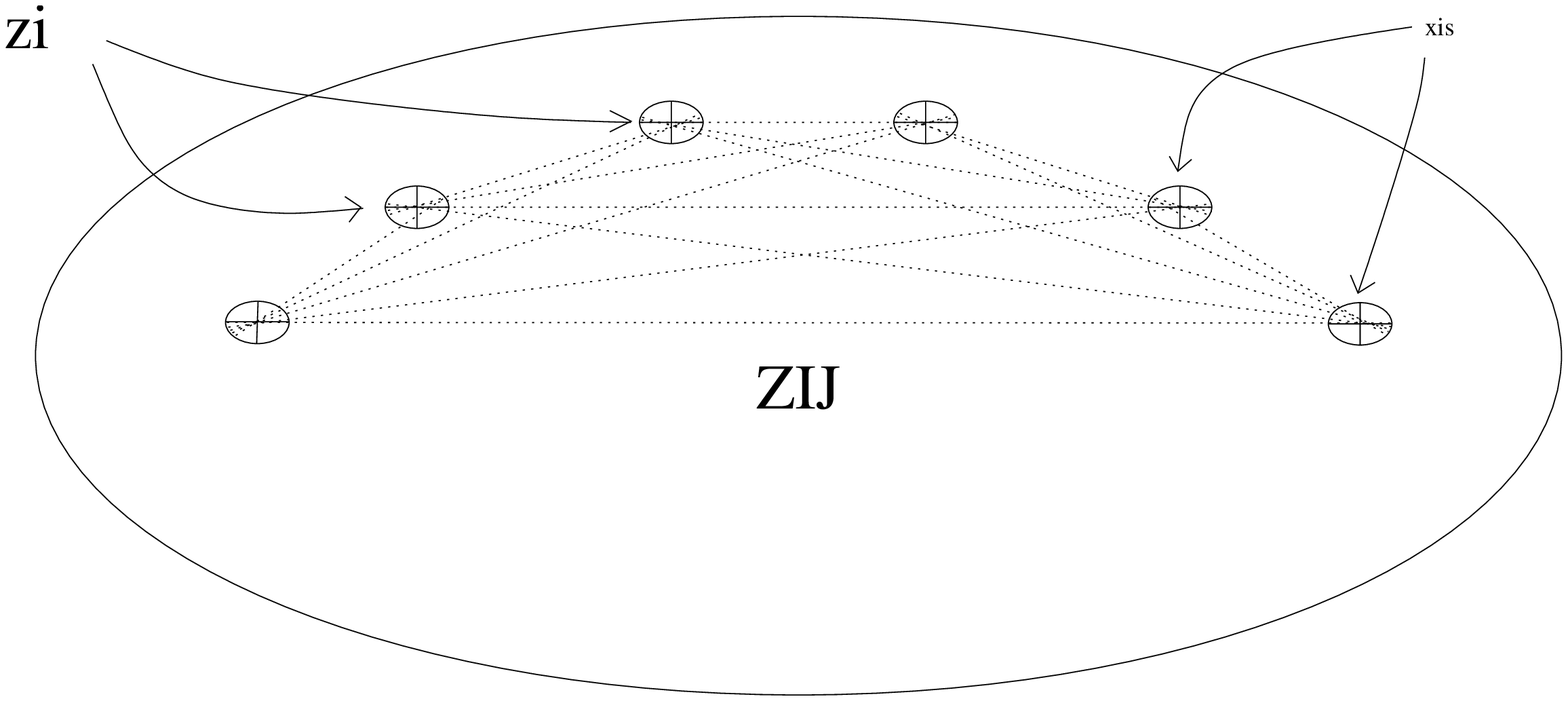} $$
Essential critical points of $d$, $n=6$.
}
\end{center}
\end{minipage}
\end{center}

The $1$-cells which correspond to the minor axis of an ellipse
will be labeled $x_i$ for $i \in \{1, \cdots, n\}$. $x_i$
can be thought of as a loop of configurations (starting at the base-point)
where one point of the configuration is stationary and the other
traverses an embedded circle which bounds a disc 
that contains only the $i$-th puncture,
in particular this disc does not contain 
the stationary point of the configuration.

The cells corresponding to the major axis of an ellipse
will be labeled $Z_i$ for $i \in \{ 1, \cdots, n \}$. $Z_i$
can be thought of as the homotopy between $y_i$ and $b^{-1}x_ibx_i$.
The cells corresponding to essential critical points $[x,y]$ where $x$
is on the $i$-th ellipse and $y$ is on the $j$-th ellipse
will be labeled $Z_{i,j}$ for $1 \leq i < j \leq n$. The attaching
map for these critical points is clearly a commutator, we choose
the attaching map to be $[bx_ib^{-1},x_j^{-1}]$ for technical convenience
in section 3.3.
\qed
\smallskip

We can reduce the CW-complex $Z$, using $n$ handle slides described by
the attaching maps for $Z_i$ for $i \in \{1, \cdots, n\}$ to get
the following CW-decomposition.

{\bf Corollary~1.}
\label{cw1a}
$\Config_2 \PDisc_n$ deformation
retracts to a subspace with CW-decomp\-osition denoted simply
by $Y$, with one $0$-cell, $n+1$ $1$-cells, labeled
$x_i : i \in \{ 1, \cdots , n \}$ and $b$, and ${n \choose 2} + n$
$2$-cells, with attaching maps given by:

$$
\begin{array}{llr}
 Y_{i,j} & [bx_ib^{-1},x_j^{-1}] & 1\leq i<j\leq n \\
 Y_i & [b,x_ibx_i] & 1\leq i\leq n
\end{array}
$$

\subsection{The dual CW-decomposition}

{\bf Lemma~3.}
\label{cw2}
The negative of the distance function $-d$
gives a CW decomposition of $\Config_2 \PDisc_n$ with $1$ $0$-cell, 
$2n+1$ $1$-cells labeled $s, a_i, b_i : i \in \{1 , \cdots, n\}$, and
${n \choose 2} + 2n$ $2$-cells, with 
attaching maps given by:

$$
\begin{array}{llr}
X_{i,j} & [a_i,b_j] & 1\leq i<j\leq n \\
X_i & s^{-1}b_i^{-1}sa_i & 1\leq i\leq n \\
X'_i & b_i ((b_{i+1} \cdots b_n s^{-1} a_1\cdots a_{i-1}) \to a_i^{-1})  & 1\leq i\leq n
\end{array}
$$

where $(x\to y) = xyx^{-1}$

{\bf Proof.}
This proof differs very little from the proof of Lemma 2.
$d$ has all the same critical points as $-d$, they only differ in their
essential critical points. There are ${n \choose 2} + 2n$ $2$-cell attachments,
$2n+1$ $1$-cell attachments and a $0$-cell. The $0$-cell is the configuration 
corresponding to the major axis of the big ellipse. The $1$-cell attachments
correspond to the minor axis of the big ellipse, labeled $s$, and 
the maxima of $-d$ on the $2$-dimensional strata where one configuration
is on the big ellipse. Label these cells $a_i$ and $b_i$ for $i \in \{1,\cdots,n\}$.

\begin{center}
\begin{minipage}{5.5cm}
\begin{center}Figure 5
{\psfrag{ai}{$a_i$'s}\psfrag{bi}{$b_i$'s}\psfrag{s}{$s$}
$$\includegraphics[width=5.5cm]{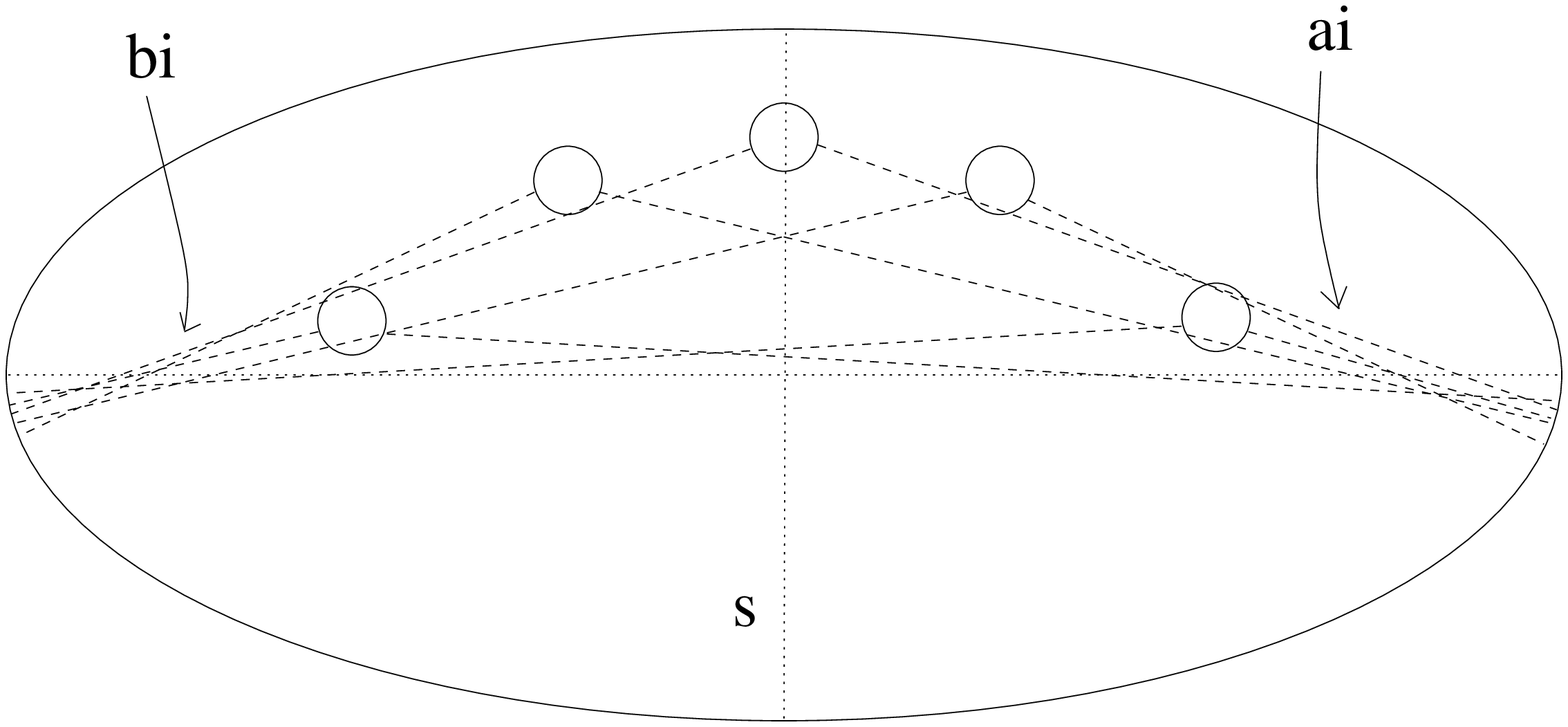} $$}
Index $0$ and $1$ critical points of $-d$, and cell
labels.
\end{center}
\end{minipage}
\hskip 3mm
\begin{minipage}{5.5cm}
\begin{center}Figure 6
{\psfrag{xp}{$X_i'$'s}\psfrag{xi}{$X_i$'s}\psfrag{XIJ}{$X_{i,j}$'s}
$$\includegraphics[width=5.5cm]{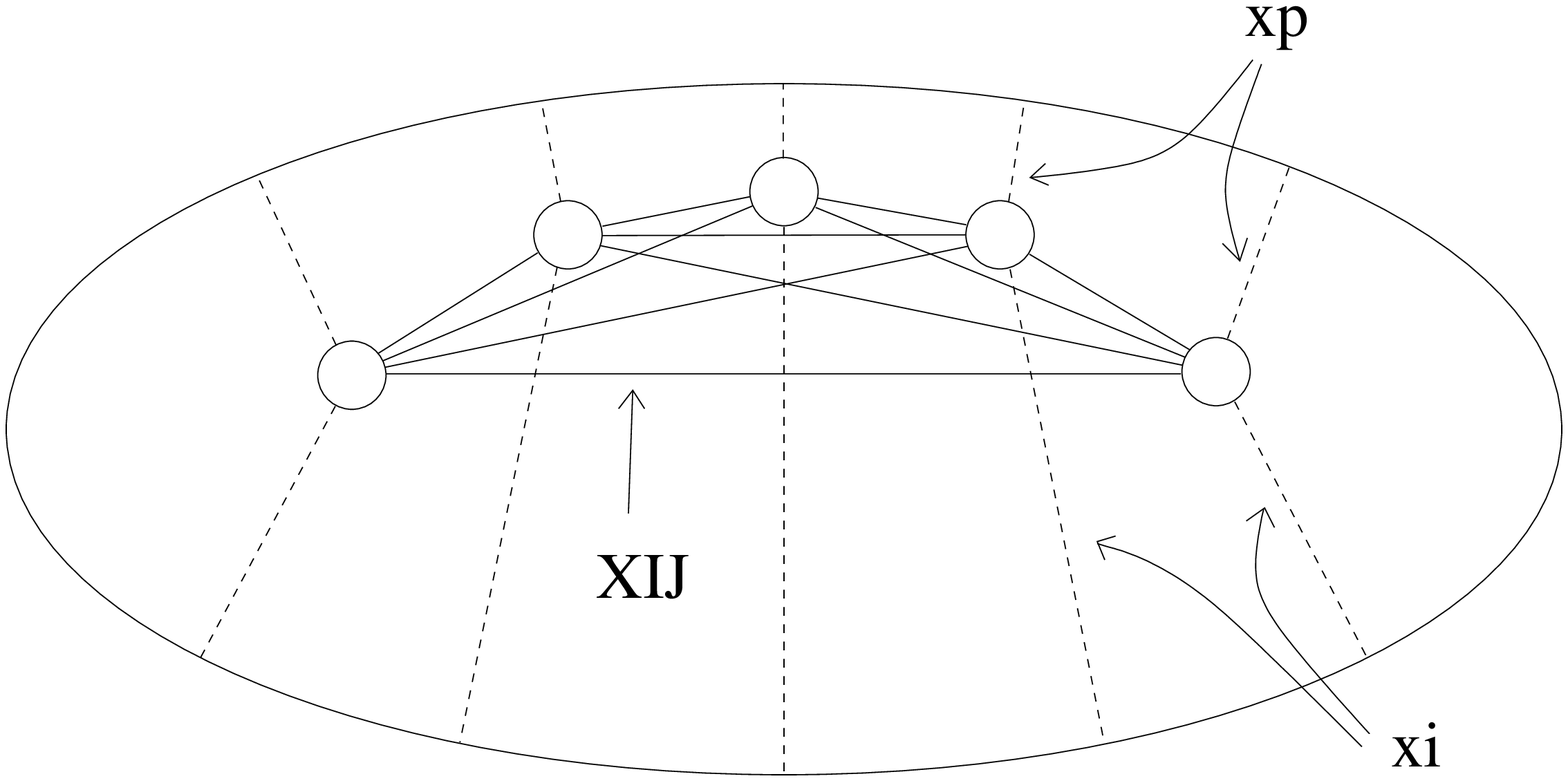} $$
Index $2$ critical points of $-d$ and cell labels.}
\end{center}
\end{minipage}
\end{center}

The $2$-cells $X_{i,j}$ are as in Corollary 1. They are described as the
attaching map for the single $2$-cell of the torus that consists of
configurations of two points, one point that lays on a circle that bounds
the $i$-th puncture and one point that lays on a circle that bounds the
$j$-th puncture.  $X_i$ and $X_i'$ express the $a_i$'s in terms of the $b_j$'s
in the two ways: conjugating by $s$ and also conjugating by its inverse.\qed

A CW-decomposition of the base space of a covering space lifts to a CW-decomposition
of the cover. Choose a base-point in $\La \Config_2 \PDisc_n$ that is above the
$0$-cell for our CW-decomposition $Y$.  We give $\La \Config_2 \PDisc_n$
a CW-structure by lifting the cells of $Y$ to $\La \Config_2 \PDisc_n$.
Introducing a mild notational ambiguity, we interpret $Y_{i,j}$ to be
the lift of $Y_{i,j}$ to $\La \Config_2 \PDisc_n$ so that the attaching map
starts at the base-point.  This is a slight abuse of notation because
$Y_{i,j}$ is also a $2$-cell in $\Config_2 \PDisc_n$.  We do the same for all the
remaining cells of $Y$.
With these conventions, $\partial Y_{i,j} = q^{-1}t^{-1}(q-1)((-x_i+tx_j)-(q-1)b)$
and $\partial Y_i = q^{-1}t^{-1}(q^{-1}t^{-1} + 1)((1-t)x_i+(q-1)b)$.

{\bf Proposition~1.}
The homology $\Ho_2 (\La \Config_2 \PDisc_n)$
contains a free rank ${n \choose 2}$ $\Lau$-module, which is
$B_n$-invariant, and spanned by
$v_{i,j} = qt(q-1)(Y_i - tY_j) + (1-t)(1+qt)Y_{i,j}$.
This submodule is the Lawrence-Krammer module and denoted
by $\La_n$. The corresponding representation of $B_n$ is called the
Lawrence-Krammer representation.

{\bf Proof.} Homotope $\sigma_i : \La \Config_2 \PDisc_n \to \La \Config_2 \PDisc_n$ to
a cellular map.  Since $\pi_2 Y = *$, we can ask how $\sigma_i$
acts on the elements $Y_j, Y_{j,k} \in \pi_2 (Y,Y^1)$, where $Y$ is the
CW-structure from Corollary 1, lifted to $\La \Config_2 \PDisc_n$.
Choose the cellular approximation so that $\sigma_i\cdot b=b$ and
$\sigma_i \cdot x_j=x_j$ unless $i=j-1$ in which case
$\sigma_i \cdot x_j=x_j^{-1}x_{j-1}x_j$ or if $i=j$ then $\sigma_i \cdot x_j=x_{j+1}$.
We can now compute $\sigma_i \cdot Y_{j,k}$ and $\sigma_i Y_j$.

$$
\sigma_i\cdot Y_{j,k} = \left\{
\begin{array}{lr}
Y_{j+1,k} & j=i, k>i+1 \\
Y_{j,k+1} & j<i, k=i \\
-qtY_{j,k} + qt(1-q)Y_k & j=i, k=i+1 \\
qY_{j-1,k} + (1-q)Y_{j,k} & j=i+1 \\
(1-q)Y_{j,k} + qY_{j,k-1} & j<i, k=i+1 \\
Y_{j,k} & \text{otherwise}
\end{array}
\right.
$$

$$
\sigma_i\cdot Y_j = \left\{
\begin{array}{lr}
Y_{j+1} & j=i \\
(1+qt)(1-q)Y_{j} + q^2Y_{j-1} + t^{-1}(1+qt)(1-t)Y_{j-1,j} & j=i+1 \\
Y_j & \text{otherwise}
\end{array}
\right.
$$

\begin{center}
\begin{minipage}{10cm}
\begin{center} Figure 7
{
\psfrag{YIP1}{$Y_{i+1}$}
\psfrag{b}[tl][tl][0.8][0]{$b$}
\psfrag{xip1}[tl][tl][0.8][0]{$x_{i+1}$}
\psfrag{xI}[tl][tl][0.8][0]{$x_i$}
\psfrag{myip}{$-Y_{i+1}$}
\psfrag{pyip}{$+Y_{i+1}$}
\psfrag{pyI}{$+Y_i$}
\psfrag{si}{$\sigma_i$}
\psfrag{pyiip1}{$+Y_{i,i+1}$}
\psfrag{myiip1}{$-Y_{i,i+1}$}
$$\includegraphics[width=10cm]{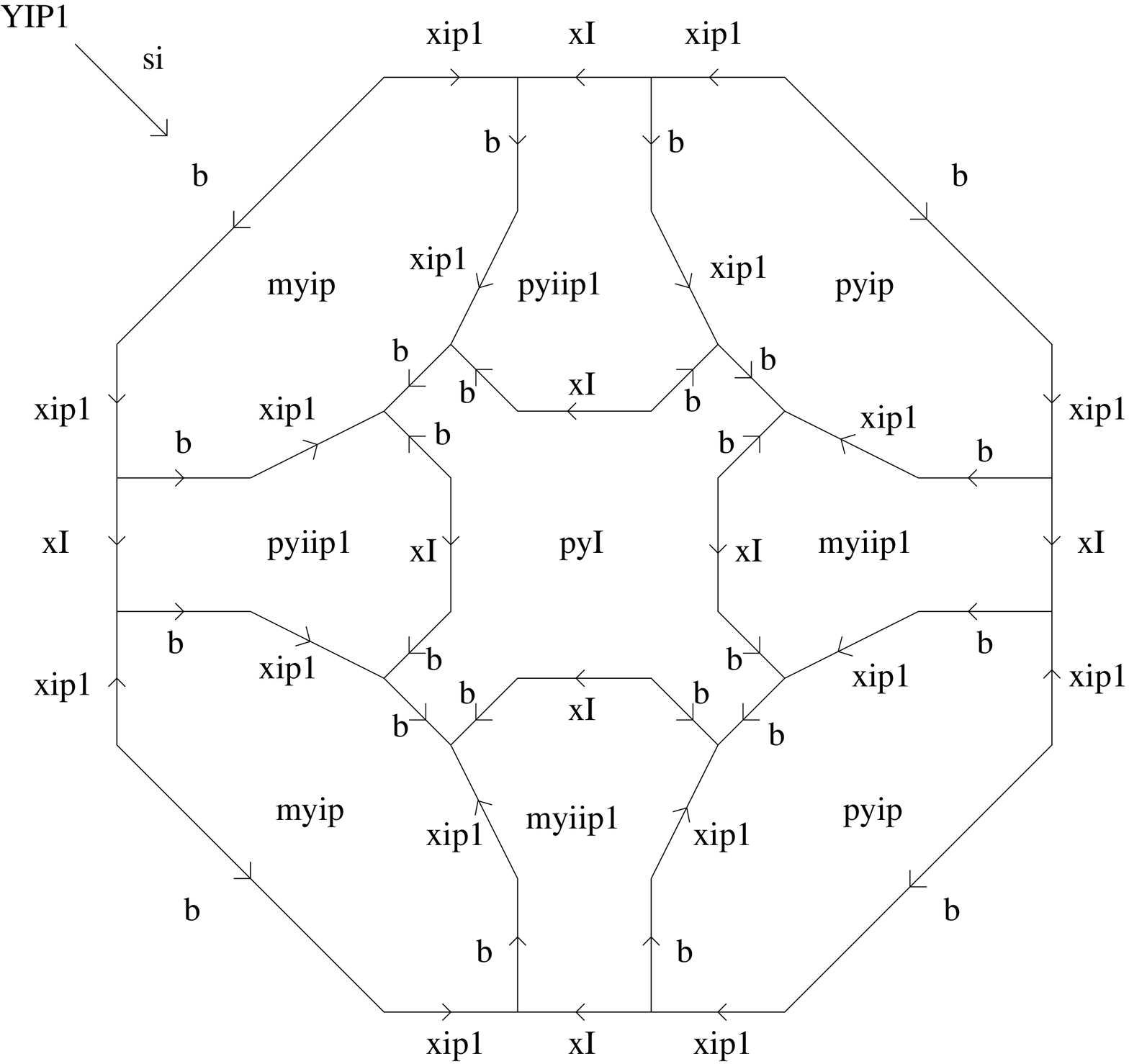}$$
Filling $\partial (\sigma_i\cdot Y_{i+1})$
}
\end{center}
\end{minipage}
\end{center}

The most involved of these computations is for $\sigma_i\cdot Y_{i+1}$ so
it will be done in some detail. The remaining computations are simpler.
$\partial Y_{i+1} = [b,x_{i+1}bx_{i+1}]$
so $\partial (\sigma_i\cdot Y_{i+1}) = [b,x_{i+1}^{-1}x_ix_{i+1}bx_{i+1}^{-1}x_ix_{i+1}]$.
Since $\pi_2 Y = *$ there is only one element of $\pi_2 (Y,Y^1)$ with this
boundary.  The problem of finding this element is much like solving
a jigsaw-puzzle, one needs to find the $2$-cells of $Y$ that fit together
so that they have boundary 
$[b,x_{i+1}^{-1}x_ix_{i+1}bx_{i+1}^{-1}x_ix_{i+1}]$.  The solution
is given in Figure 7.

The above computation proves more than the fact that subspace spanned by the
$v_{i,j}$'s is invariant. It gives us the matrices
for the Lawrence-Krammer representation. They are:
$$
\sigma_i\cdot v_{j,k} = \left\{
\begin{array}{lr}
v_{j,k} & i\notin \{j-1,j,k-1,k\}, \\
qv_{i,k} + (q^2-q)v_{i,j} + (1-q)v_{j,k} & i=j-1 \\
v_{j+1,k} & i=j\neq k-1, \\
qv_{j,i} + (1-q)v_{j,k} - (q^2-q)tv_{i,k} & i=k-1\neq j,\\
v_{j,k+1} & i=k,\\
-tq^2v_{j,k} & i=j=k-1.
\end{array}
\right.
$$
\qed

The above computation gives an independent verification that the matrices in
\cite{BB} are correct.

Even though the $v_{i,j}$'s form a basis for
$\Ho_2(\La \Config_2 \PDisc_n) \otimes_{\Lau} \F$ where $\F$ is a field
containing $\Lau$, the $v_{i,j}$'s do not
span $\Ho_2(\La \Config_2 \PDisc_n)$.  For example, the homology
classes given by Bigelow in \cite{Bi3} are not in the $\Lau$-span
of the $v_{i,j}$'s. This is similar to, although much less complete,
than the results of Paoluzzi and Paris \cite{Pl},
 where they prove that $\Ho_2(\La \Config_2 \PDisc_n)$
is a free rank $n \choose 2$ $\Lau$-module.

\section{On the image of the Lawrence-Krammer representation}
\label{ILK}

It is not uncommon for a general manifold to have a singular intersection product
pairing.
For example, with the cylinder $S^1 \times \I$ the intersection product is
zero. There is a Poincar\'e duality theorem for abelian covers of
compact manifolds (See for example \cite{Kaw} Appendix E), but like Poincar\'e
duality for compact manifolds, Poincar\'e duality does not directly give information
about the intersection product provided the manifolds have non-empty
boundary.

To compute the intersection product, one could take the CW-decomposition
given in \cite{Bi1} and notice that all homology classes are realizable
by compact, genus two surfaces.  General position is sufficient to compute
the pairings.  Unfortunately, Bigelow's $v_{i,j}$'s are not transverse.
Instead, we compute the intersection product on $\Ho_2 (\La \Config_2 \PDisc_n)$
using the two transverse CW-decompositions, $X$ and $Y$.

We lift the cells of $X$ to $\La \Config_2 \PDisc_n$ as we did for $Y$,
that is, fix some lift of $X^0$ to $\La \Config_2 \PDisc_n$ and call it
the base-point. For every cell $X_{i,j}$, $X_i$ and $X_i'$ of $X$, we use
the same notation to denote the lift of the cell in $\La \Config_2 \PDisc_n$
whose attaching map starts at $X^0$.

{\bf Theorem~2.}
The intersection product $\langle v_{i,j}, v_{k,l}\rangle$ is given by the
formula
$$
-(1-t)(1+qt)(q-1)^2t^{-2}q^{-3}
\left\{
\begin{array}{lr}
-q^2t^2(q-1) & i=k<j<l \text{ or } i<k<j=l \\
-(q-1) & k=i<l<j \text{ or } k<i<j=l \\
t(q-1) & i<j=k<l \\
q^2t(q-1) & k<l=i<j \\
-t(q-1)^2(1+qt) & i<k<j<l \\
(q-1)^2(1+qt) & k<i<l<j \\
(1-qt)(1+q^2t) & k=i, j=l \\
0 & \text{otherwise} \\
\end{array}
\right.
$$

{\bf Proof.}
To see that the above formula is correct, notice that the two 
CW-decomp\-osi\-tions
given for $\Config_2 \PDisc_n$ in Corollary 1 and Lemma 3 are transverse.
In fact, the only cells
of $X$ that intersect cells of $Y$ are the ${n \choose 2}$ pairs,
$X_{i,j} \cap Y_{i,j}$, which intersect in precisely four points (before lifting
to the cover). This
can be seen easily because the $2$-cells $X_i$ and $X_i'$ are all contained
in the $3$-dimensional stratum of $\Config_2 \PDisc_n$, and $X_{i,j}$ is
disjoint from $Y_{k,l}$ unless $k=i$ and $j=l$ as in Figure 8. The disjointness
observation comes from the Morse Theory of \cite{Ha2} -- any cell in a
CW-decomposition for a Morse function $f$ can be realized in a very simple
way: a Morse function near a critical point $z$ has a local coordinate
system where $f(x_1,\cdots,x_n)=-(x_1^2+\cdots+x_k^2)+(x_{k+1}^2+\cdots+x_n^2)+f(z)$,
where $k$ is the index of the critical point.  Let $D$ be a compact $k$-dimensional
disc that corresponds to a neighborhood of $0$ inside the subspace $x_{k+1}=\cdots=x_n=0$
in the above coordinate system.  Then the cell corresponding to the critical
point $z$ consists of $D$ union the forward orbit of $\partial D$ under
the flow of the negative gradient of the Morse function.

\begin{center}
\begin{minipage}{9cm}
\begin{center} Figure 8
{
\psfrag{xi}{$2$}
\psfrag{xj}{$5$}
\psfrag{xij}{$X_{2,5}$}
\psfrag{yk}{$1$}
\psfrag{yl}{$4$}
\psfrag{ykl}{$Y_{1,4}$}
$$\includegraphics[width=9cm]{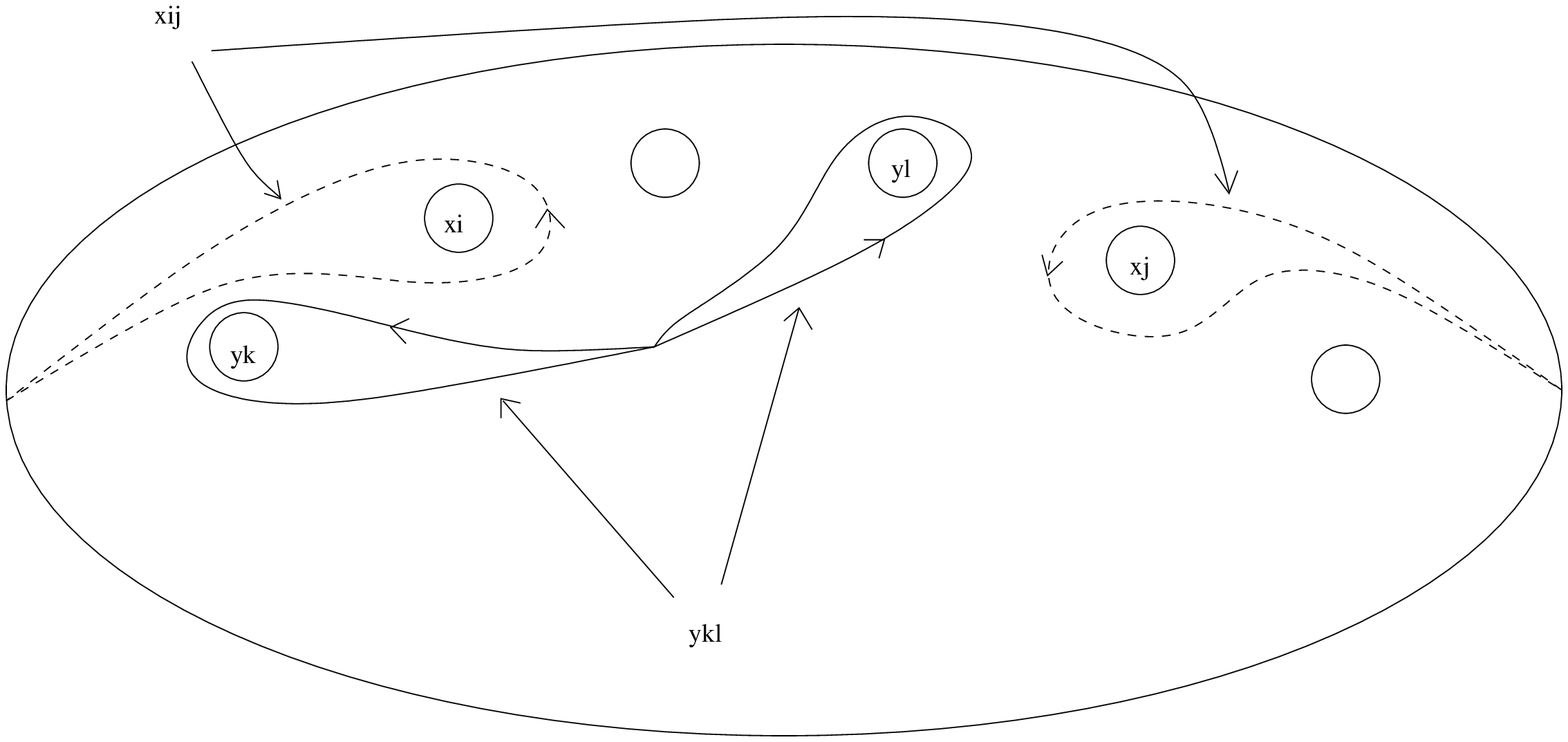} $$
}
$2$-cells $Y_{1,4}$ and $X_{2,5}$, $n=6$.
\end{center}
\end{minipage}
\end{center}

An easy computation gives
$\langle X_{i,j}, Y_{i,j} \rangle = q^{n+3-(i+j)}(q-1)^2$.
To compute $\langle v_{i,j}, v_{k,l}\rangle$ we represent $v_{i,j}$
in the $X$ cellular homology and use the previous formula
for $\langle X_{i,j},Y_{k,l} \rangle$. Since both $X$ and $Y$ are
CW-decompositions of deformation retractions of $\Config_2 \PDisc_n$,
they are canonically homotopy equivalent.  The homotopy equivalence
$Y \to X$ can be homotoped to a cellular map where
$x_i \to (b_n \cdots b_{i+1}) b_i^{-1} (b_n \cdots b_{i+1})^{-1}$,
and $b \to s(b_n \cdots b_1)^{-1}$. Using Lemma 4, we can compute
the homotopy equivalence on the $2$-cells of $Y$ as in the proof of
Proposition 1.

$$ Y_{i,j} = q^{-4} q^{i+j-n} X_{i,j} -
   \sum_{\substack{0<a<i \\ j<b<n+1}} q^{-5}(q-1)^2 q^{a+b-n} X_{a,b} +
   \sum_a l_a X_a +$$
$$ \sum_{j<b<n+1} q^{-5}(q-1) q^{i+b-n} X_{i,b} +
   \sum_{0<a<i} -q^{-4}(q-1) q^{a+j-n} X_{a,j} $$

See Figure 9 for details on the above computation for $Y_{i,j}$.

\begin{center}
\begin{minipage}{13cm}
\begin{center} Figure 9
{
\psfrag{bi}[tl][tl][0.6][41]{$b_i$}
\psfrag{bil}[tl][tl][0.6][-49]{$b_i$}
\psfrag{ai}[tl][tl][0.6][-49]{$a_i$}
\psfrag{YIJ}{$Y_{i,j}$}
\psfrag{s}[tl][tl][0.6][41]{$s$}
\psfrag{aim1}[tl][tl][0.6][-49]{$a_{i-1}\cdots a_1$}
\psfrag{bnj40}[tl][tl][0.6][41]{$b_{n}\cdots b_{j+1}$}
\psfrag{bimb1}[tl][tl][0.6][-49]{$b_{i-1}\cdots b_1$}
\psfrag{xij}[tl][tl][0.8][-49]{$+X_{i,j}$}
\psfrag{px*j}[tl][tl][0.8][-49]{$+X_{*,j}$}
\psfrag{mx*j}[tl][tl][0.8][-49]{$-X_{*,j}$}
\psfrag{pxi*}[tl][tl][0.8][-49]{$+X_{i,*}$}
\psfrag{mxi*}[tl][tl][0.8][-49]{$-X_{i,*}$}
\psfrag{px**}[tl][tl][0.8][-49]{$+X_{*,*}$}
\psfrag{mx**}[tl][tl][0.8][-49]{$-X_{*,*}$}
\psfrag{xi}[tl][tl][0.8][-49]{$X_i$}
\psfrag{x*}[tl][tl][0.8][-49]{$X_*$}
\psfrag{bnbjp1}[tl][tl][0.7][17]{$b_{n}\cdots b_{j+1}$}
\psfrag{bnbjp8}[tl][tl][0.7][57]{$b_{n}\cdots b_{j+1}$}
\psfrag{bnbjp10}[tl][tl][0.7][20]{$b_{n}\cdots b_{j+1}$}
\psfrag{bnbjp18}[tl][tl][0.7][60]{$b_{n}\cdots b_{j+1}$}
\psfrag{bnb12}[tl][tl][0.7][-15]{$b_{n}\cdots b_{1}$}
\psfrag{bnb16}[tl][tl][0.7][-85]{$b_{n}\cdots b_{1}$}
\psfrag{bnb112}[tl][tl][0.7][-20]{$b_{n}\cdots b_{1}$}
\psfrag{bnb116}[tl][tl][0.7][-90]{$b_{n}\cdots b_{1}$}
\psfrag{bnbip3}[tl][tl][0.7][-30]{$b_{n}\cdots b_{i+1}$}
\psfrag{bnbip5}[tl][tl][0.7][-68]{$b_{n}\cdots b_{i+1}$}
\psfrag{bnbip13}[tl][tl][0.7][-30]{$b_{n}\cdots b_{i+1}$}
\psfrag{bnbip15}[tl][tl][0.7][-68]{$b_{n}\cdots b_{i+1}$}
\psfrag{bi4}[tl][tl][0.7][-45]{$b_i$}
\psfrag{bi9}[tl][tl][0.7][45]{$b_i$}
\psfrag{bi14}[tl][tl][0.7][-45]{$b_i$}
\psfrag{bi19}[tl][tl][0.7][45]{$b_i$}
\psfrag{s20}[tl][tl][0.7][10]{$s$}
\psfrag{s7}[tl][tl][0.7][80]{$s$}
\psfrag{s11}[tl][tl][0.7][5]{$s$}
\psfrag{s17}[tl][tl][0.7][80]{$s$}
$$\includegraphics[width=13cm]{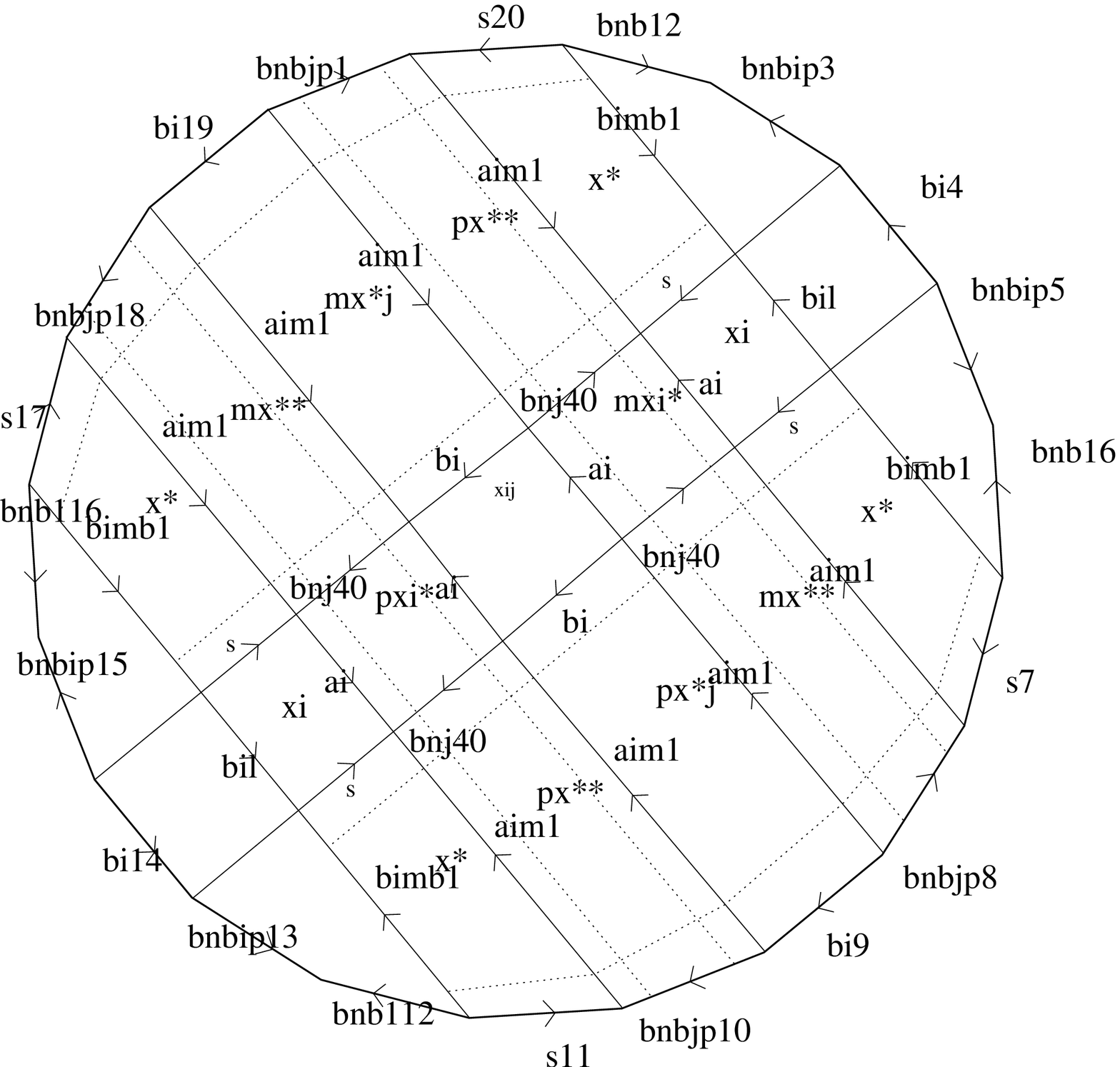} $$
}
Filling $\partial Y_{i,j} = [s(b_n\cdots b_1)^{-1}\to (b_n\cdots b_{i+1} \to b_i^{-1}),b_n\cdots b_{j+1} \to b_j]$
\end{center}
\end{minipage}
\end{center}

$$ Y_i = -\sum_{\substack{0<a<i \\ i<b<n+1}} q^{-6}t^{-1}(q-1)(1+qt) q^{a+b-n} X_{a,b} +
        \sum_a p_a X_a + $$
$$ \sum_{0<a<i} q^{-4} q^{a+i-n} X_{a,i} -
        \sum_{i<b<n+1} q^{-6}t^{-1} q^{i+b-n} X_{i,b}$$

See Figure 10 for details on the computation for $Y_i$.

There is no need to compute the coefficients $l_a$ and
$p_a$ as they do not contribute to the intersection product.

\begin{center}
\begin{minipage}{13cm}
\begin{center} Figure 10
{
\psfrag{Yi}{$Y_i$}
\psfrag{xipr}[tl][tl][0.9][21]{$X_i'$}
\psfrag{x*}[tl][tl][0.8][0]{$X_*$}
\psfrag{xi}[tl][tl][0.8][0]{$X_i$}
\psfrag{mx*i}[tl][tl][0.8][0]{$-X_{*,i}$}
\psfrag{mx*d*}[tl][tl][0.8][0]{$-X_{*,*}$}
\psfrag{px*d*}[tl][tl][0.8][0]{$+X_{*,*}$}
\psfrag{px*di}[tl][tl][0.8][0]{$+X_{*,i}$}
\psfrag{bi}[tl][tl][0.6][0]{$b_i$}
\psfrag{ai}[tl][tl][0.6][0]{$a_i$}
\psfrag{bntb1}[tl][tl][0.8][0]{$b_n \cdots b_1$}
\psfrag{bnbip}[tl][tl][0.8][0]{$b_n \cdots b_{i+1}$}
\psfrag{aim1}[tl][tl][0.6][-80]{$a_{i-1} \cdots a_1$}
\psfrag{aim2}[tl][tl][0.6][-58]{$a_{i-1} \cdots a_1$}
\psfrag{bnb1}[tl][tl][0.6][-45]{$b_{n} \cdots b_{i+1}$}
\psfrag{bnb2}[tl][tl][0.6][-85]{$b_{n} \cdots b_{i+1}$}
\psfrag{bnbim10}[tl][tl][0.6][5]{$b_n \cdots b_{i-1}$}
\psfrag{bnbim45}[tl][tl][0.6][40]{$b_n \cdots b_{i-1}$}
\psfrag{anaim20}[tl][tl][0.6][15]{$a_n \cdots a_{i-1}$}
\psfrag{anaim30}[tl][tl][0.6][30]{$a_n \cdots a_{i-1}$}
\psfrag{anaim25}[tl][tl][0.6][21]{$a_n \cdots a_{i-1}$}
\psfrag{bnbi23}[tl][tl][0.6][23]{$b_n \cdots b_{i+1}$}
\psfrag{bim2}[tl][tl][0.6][-55]{$b_{i-1} \cdots b_1$}
\psfrag{bim1}[tl][tl][0.6][-80]{$b_{i-1} \cdots b_1$}
\psfrag{bnbip2}[tl][tl][0.7][2]{$b_n \cdots b_{i+1}$}
\psfrag{bnbip5}[tl][tl][0.7][-50]{$b_n \cdots b_{i+1}$}
\psfrag{bnbip7}[tl][tl][0.7][-85]{$b_n \cdots b_{i+1}$}
\psfrag{bnbip10}[tl][tl][0.7][40]{$b_n \cdots b_{i+1}$}
\psfrag{bnbip12}[tl][tl][0.7][5]{$b_n \cdots b_{i+1}$}
\psfrag{bnbip15}[tl][tl][0.7][-50]{$b_n \cdots b_{i+1}$}
\psfrag{bnbip17}[tl][tl][0.7][-85]{$b_n \cdots b_{i+1}$}
\psfrag{bnbip20}[tl][tl][0.7][40]{$b_n \cdots b_{i+1}$}
\psfrag{bi1}[tl][tl][0.7][25]{$b_i$}
\psfrag{bi6}[tl][tl][0.7][-65]{$b_i$}
\psfrag{bi11}[tl][tl][0.7][30]{$b_i$}
\psfrag{bi16}[tl][tl][0.7][-65]{$b_i$}
\psfrag{bntb13}[tl][tl][0.7][-15]{$b_n \cdots b_1$}
\psfrag{bntb19}[tl][tl][0.7][55]{$b_n \cdots b_1$}
\psfrag{bntb114}[tl][tl][0.7][-30]{$b_n \cdots b_1$}
\psfrag{bntb118}[tl][tl][0.7][-105]{$b_n \cdots b_1$}
\psfrag{s4}[tl][tl][0.7][-30]{$s$}
\psfrag{s8}[tl][tl][0.7][75]{$s$}
\psfrag{s13}[tl][tl][0.7][-10]{$s$}
\psfrag{s19}[tl][tl][0.7][70]{$s$}
\psfrag{sa}[tl][tl][0.7][21]{$s$}
\psfrag{sb}[tl][tl][0.7][21]{$s$}
\psfrag{sc}[tl][tl][0.7][-30]{$s$}
\psfrag{sd}[tl][tl][0.7][70]{$s$}
\psfrag{mx*d*2}[tl][tl][0.8][21]{$-X_{*,*}$}
\psfrag{px*d*2}[tl][tl][0.8][21]{$+X_{*,*}$}
\psfrag{px*di2}[tl][tl][0.8][21]{$+X_{*,i}$}
\psfrag{x*2}[tl][tl][0.8][15]{$X_*$}
\psfrag{x*3}[tl][tl][0.8][26]{$X_*$}
\psfrag{mx*d*tr}[tl][tl][0.8][15]{$-X_{*,*}$}
\psfrag{px*d*br}[tl][tl][0.8][26]{$+X_{*,*}$}
\psfrag{mx*imr}[tl][tl][0.8][21]{$-X_{*,i}$}
\psfrag{x*tl}[tl][tl][0.8][15]{$X_*$}
\psfrag{x*bl}[tl][tl][0.8][26]{$X_*$}
\psfrag{ximl}[tl][tl][0.8][21]{$X_i$}
$$\includegraphics[width=13cm]{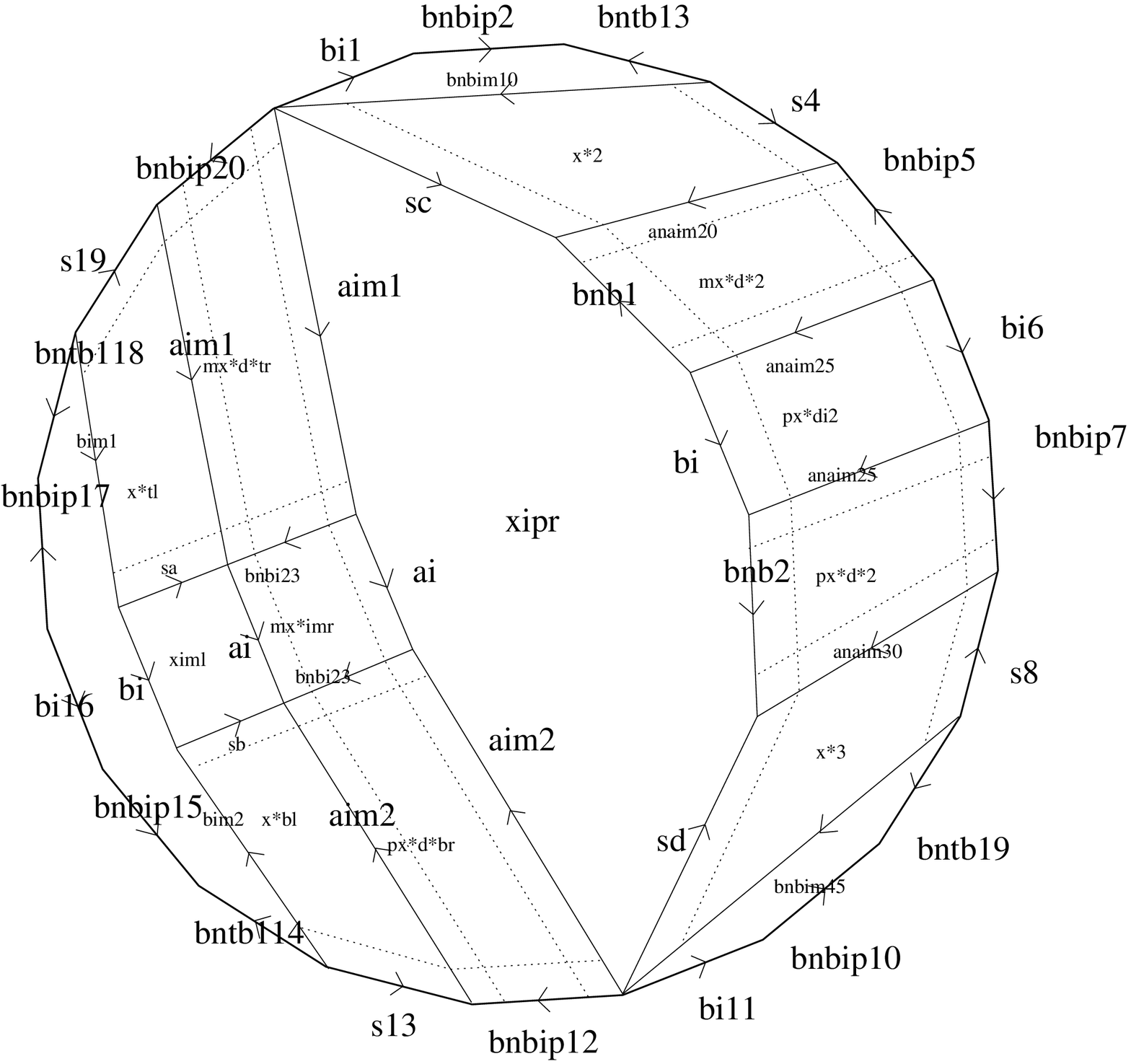} $$
}
$\partial Y_{i} = [s(b_n\cdots b_1)^{-1},(b_n\cdots b_{i+1} \to b_i^{-1})s(b_n\cdots b_1)^{-1}(b_n \cdots b_{i+1} \to b_i^{-1})]$
\end{center}
\end{minipage}
\end{center}
 \qed

{\bf Proposition~2.}
The intersection product is non-singular,
or equivalently, the dual map $v \to \langle \cdot, v \rangle$ is injective.

{\bf Proof.}
We need to prove that the ${n \choose 2} \times {n \choose 2}$ matrix of
coefficients $\langle v_{i,j},v_{k,l}\rangle$ has rank $n \choose 2$, thus
it suffices to show that the determinant of this matrix is non-zero.
Since $\Lau$ is an integral domain (moreover it is a 
unique factorization domain
\cite{Kaw}), we really only need to show that
the determinant of the matrix $M_{(i,j),(k,l)}=\frac{\langle v_{i,j},v_{k,l}\rangle}{c_{p,q}}$ is
non-zero where $c_{p,q}=-(1-t)(1+qt)(q-1)^{2}t^{-2}q^{-3}$. 
Notice the term in the determinant of
$M$ corresponding to the diagonal entries is precisely 
$\left( (1-qt)(1+q^2t) \right)^{n \choose 2}$.
Similarly notice that every other term in the determinant of $M$ is divisible by
$(q-1)$.  Since neither $(1-qt)$ or $(1+q^2t)$ are divisible by $(q-1)$, the determinant
must be non-zero.
\qed

{\bf Theorem~3.}
\label{theoremunitary}
For appropriate choices of $q$ and $t \in \C$
the intersection product is a negative-definite Hermitian form.

{\bf Proof.}
Take an arbitrary $v \in \La_n$ and compute
$\langle v,v \rangle$.  Let $v = \sum_{i,j} \lambda_{i,j}v_{i,j}$
for $\lambda_{i,j} \in \C$, and notice that for $|q|=|t|=1$,
$-(1-t)(1-q^2t^2)(q-1)^2(1+q^2t)t^{-2}q^{-3} \in \Real$.
By Theorem 2

$$\frac{\langle v,v \rangle}{-(1-t)(1-q^2t^2)(q-1)^2(1+q^2t)t^{-2}q^{-3}}
 = \sum_{i,j} \lambda_{i,j}\bar\lambda_{i,j} + 2\text{Re}\left( \frac{q-1}{(1-qt)(1+q^2t)}k\right)$$
where
$$k=
\sum_{\substack{a=c, b>d \\ or \\ b=d, c<a}}\lambda_{a,b}\bar\lambda_{c,d} +
\sum_{a=d} q^2t\lambda_{a,b}\bar\lambda_{c,d}+
\sum_{c<a<d<b} (q-1)(1+qt)\lambda_{a,b}\bar\lambda_{c,d}
$$

\noindent therefore, for $|q-1|<\frac{1}{2n^4+6n^3}$ and $|t-i|<\frac{1}{2n^4+6n^3}$
 the Lawrence-Krammer representation is definite. Note that
$-(1-t)(1-q^2t^2)(q-1)^2(1+q^2t)t^{-2}q^{-3}$ is negative
in this case.\qed

\section{Conjugacy in the image}

Since the Lawrence-Krammer representation is faithful, one may ask if
it gives insight into the conjugacy problem for braid groups.
One way to approach this would be via canonical forms of matrices.
Given a braid $f \in B_n$ let $f_*$ denote the action of $f$ on $\La_n$.
Unitary matrices can be diagonalized, and diagonal matrices are conjugate
if an only if they have the same characteristic polynomial.
If two matrices in the image of the
Lawrence-Krammer representation are conjugate in $U_{n \choose 2}$, are
they conjugate by a matrix in the image of the Lawrence-Krammer
representation?  If the answer is yes, this would be an exceptionally
fast solution to the conjugacy problem in braid groups. It turns out the
answer is no, and this will be proved in Corollary 2.

The fact that the Lawrence-Krammer representation is unitary, together
with its simple topological definition allows the proof of certain
symmetry relations among the eigenvalues of $f_*$ and the eigenvalues
of the matrices of related braids.  These symmetries,
together with the existence of non-invertible knots \cite{Tr} will be used
to show that the characteristic polynomial does not separate conjugacy
classes.

Given a braid $f \in B_n$ there is an associated braid
$cfc$, where $c : \PDisc_n \to \PDisc_n$
is any orientation-reversing diffeomorphism
of the punctured disc that fixes the $n$ puncture points.
$c$ can be chosen to have order $2$. If one thinks of
$\PDisc_n$ as the unit disc in the complex plane
with puncture points along the real axis, then $c$ can be
taken to be complex conjugation.
The map from $B_n \to B_n$ given
by $f \to cfc$ is the only outer automorphism of the braid groups
$B_n$ \cite{DG}.

{\bf Proposition~3.}
\label{conjf} The matrices $f_*$ and $(cf^{-1}c)_*$ are conjugate
in $U_{n \choose 2}$.

{\bf Proof.}
Notice that since $c$ can be realized as complex conjugation on
the punctured disc, $c$ defines an involution of $\Config_2 \PDisc_n$ which
lifts to an involution of $\La \Config_2 \PDisc_n$.
Notice that the induced map $c_*$ on $\La_n$ is not linear, in
fact $c_*(t^aq^bv) = t^{-a}q^{-b}c_*(v)$. Therefore,
$(cfc)_* = c_*f_*c_* = (c_*\tilde{I})\circ(\tilde{I} f_* \tilde{I})\circ(\tilde{I} c_*)$
is a composite of three $\Lau$-linear maps,
where $\tilde{I} : \La_n \to \La_n$ is the unique $\Zed$-linear map
such that $\tilde{I}(v_{i,j}) = v_{i,j}$ and
$\tilde{I}(t^aq^bv)=t^{-a}q^{-b}\tilde{I}(v)$.
Notice $\tilde{I} f_* \tilde{I} = \overline{f_*}$.
This proves that $\overline{f_*}$ is conjugate to
$(cfc)_*$. Since the Lawrence-Krammer representation is unitary,
$\overline{f_*}$ is conjugate
to $(f_*^{-1})^\tau$, where $\tau$ denotes the transpose operation.
The matrix that conjugates the one to the other is the matrix of
products $\langle v_{i,j},v_{k,l}\rangle$. Therefore, $(cf^{-1}c)_*$ is
conjugate to $f_*^{\tau}$, but $f_*^{\tau}$ and $f_*$ have the same
characteristic polynomials and are therefore conjugate.\qed

If we think of the closed braids associated to the four
braids $f$, $cfc$, $cf^{-1}c$ and $f^{-1}$, then the links associated to
$f$ and $f^{-1}$ are mirror reflections of each other, and the links
associated to $cfc$ and $f$ are also mirror reflections of each other.
The mirror reflection $f \to f^{-1}$ changes the orientation of the
knot, while $f \to cfc$ preserves the orientation. The oriented knots associated
to $f$ and $cf^{-1}c$ are inverses of each other.

\begin{center}
\begin{minipage}{5cm}
\begin{center} Figure 11
{
$$\includegraphics[width=5cm]{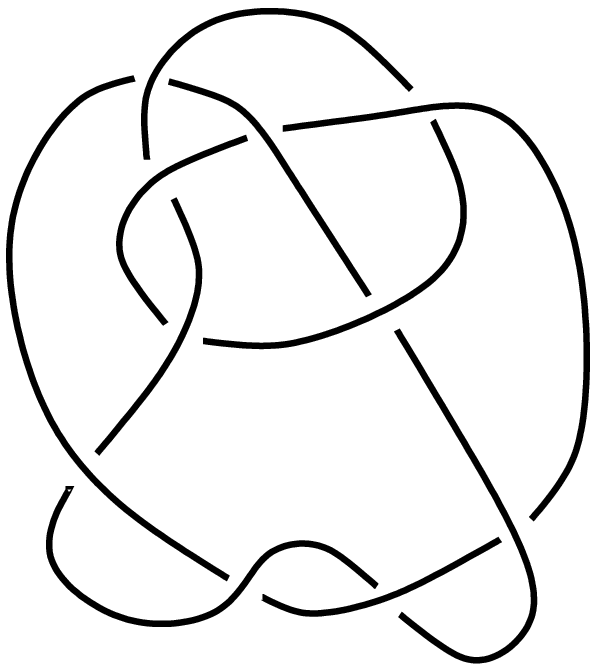} $$
}
The knot $10_{82}$
\end{center}
\end{minipage}
\end{center}

{\bf Corollary~2.}
There exists matrices in the image of the Lawrence-Krammer
representation that are conjugate, yet the braids they are associated to
are not conjugate.

Trotter has shown that non-invertible knots exist \cite{Tr}. Trotter's
example is
a rather complicated pretzel knot.  With the advent of sophisticated
computer algorithms such as Jeff Week's Snappea,
 simpler non-invertible knots have
been found \cite{HW}. For example the hyperbolic knot that is denoted
$10_{82}$ in Rolfsen's knot tables.

More generally, provided all of the Lawrence
representations \cite{La} \cite{Bi3} are unitary the above
proof would go through to prove that $f_*$ is conjugate
to $(cf^{-1}c)_*$ for all braids $f$
and all Lawrence representations. As mentioned earlier, implicit in the
work of Long \cite{Long}, the Lawrence representations all preserve
a non-singular sesquilinear form and therefore the characteristic
polynomials of $f_*$ and $(cf^{-1}c)_*$ are the same for all
Lawrence representations by the proof of Proposition 3.

\end{document}